\documentclass[11pt,reqno]{amsart}
\usepackage[T1]{fontenc}
\usepackage{lmodern}
\usepackage{amsmath,amssymb,amsthm,mathtools}
\usepackage[margin=1in]{geometry}
\usepackage{microtype}
\usepackage{xcolor}
\definecolor{linkblue}{RGB}{25,66,110}
\usepackage{hyperref}
\hypersetup{colorlinks=true,linkcolor=linkblue,citecolor=linkblue,urlcolor=linkblue,
  pdftitle={Asymptotic bounds for the Erdos-Ginzburg-Ziv constant in fixed rank},
  pdfauthor={Dongchun Han and Hanbin Zhang}}
\numberwithin{equation}{section}
\newtheorem{theorem}{Theorem}[section]
\newtheorem{proposition}[theorem]{Proposition}
\newtheorem{lemma}[theorem]{Lemma}

\newtheorem*{conj*}{Conjecture}
\newtheorem{corollary}[theorem]{Corollary}
\theoremstyle{definition}
\newtheorem{definition}[theorem]{Definition}
\theoremstyle{remark}

\newcommand{\N}{\mathbb N}
\newcommand{\Z}{\mathbb Z}
\newcommand{\Q}{\mathbb Q}
\newcommand{\R}{\mathbb R}
\newcommand{\F}{\mathbb F}
\newcommand{\T}{\mathbb T}

\newcommand{\cS}{\mathcal S}
\DeclareMathOperator{\supp}{supp}
\DeclareMathOperator{\conv}{conv}

\DeclareMathOperator{\lcm}{lcm}
\DeclareMathOperator{\im}{im}
\DeclareMathOperator{\vol}{vol}

\newcommand{\ssum}{\Sigma^{0}}

\title[Asymptotic EGZ bounds]{On the Erd\H{o}s--Ginzburg--Ziv problem for $\mathbb Z_n^3$}

\author{Dongchun Han}
\address{School of Mathematics, Southwest Jiaotong University, Chengdu 610000, P.R. China}
\email{handongchun@swjtu.edu.cn}

\author{Hanbin Zhang}
\address{School of Mathematics (Zhuhai), Sun Yat-sen University, Zhuhai 519082, Guangdong, P.R. China}
\email{zhanghb68@mail.sysu.edu.cn}

\date{}
\subjclass[2020]{11B30, 11B50, 20K01}
\keywords{Zero-sum sequence, Erd\H{o}s--Ginzburg--Ziv constant, finite abelian group, restricted sumset}

\begin{document}

\begin{abstract}
The Erd\H{o}s--Ginzburg--Ziv constant
$\mathsf s(\mathbb Z_n^r)$ is the least positive integer $t$
such that every sequence of $t$ elements of
$\mathbb Z_n^r$ contains $n$ terms whose sum is zero. In this paper, we prove that
$$\mathsf s(\Z_n^3)=9n+o(n)$$
for odd $n$. For every fixed $r\ge 4$, we also show that $\mathsf s(\Z_n^r)\le\bigl(\binom{2r-1}{r}+1\bigr)n+o_r(n).$ Our approach combines Alon--Dubiner-type spectral expansion, Zakharov's convex-geometric framework, and restricted-sum ideas developed by Han and Ren.
\end{abstract}

\maketitle

\section{Introduction}

The \emph{Erd\H{o}s--Ginzburg--Ziv constant} of $\Z_n^r$, denoted by $\mathsf s(\Z_n^r)$, is the least integer $t$ such that every sequence of $t$ elements of $\Z_n^r$, with repetitions allowed, contains a zero-sum subsequence of length $n$. The classical theorem of Erd\H{o}s, Ginzburg and Ziv~\cite{EGZ} shows that  $\mathsf s(\Z_n)=2n-1$. Later, Reiher \cite{Reiher} proved that $\mathsf s(\Z_n^2)=4n-3$. For higher ranks, the precise values of $\mathsf s(\Z_n^r)$ are largely unknown. For every fixed $r\ge2$, Alon and Dubiner~\cite{AD} proved that $\mathsf s(\Z_n^r)\le (cr\log r)^r n$, where $c$ is an absolute constant. More recently, the authors~\cite[Theorem~1.1]{HZ} established the following uniform bound:
\begin{equation}\label{eq:uniform-background}
 \mathsf s(\Z_n^r)\le\binom{2r}{r}(n-1)+1
 \qquad(n,r\ge1).
\end{equation}
On the other hand, Harborth~\cite{Harborth} proved that $\mathsf s(\Z_n^r)\ge 2^r(n-1)+1$. In rank three, Elsholtz~\cite{Elsholtz} obtained the stronger lower bound $\mathsf s(\Z_n^3)\ge9n-8$ for every odd $n\ge3$. This leads to the following natural conjecture.

\begin{conj*}
For every odd positive integer $n$, we have $\mathsf s(\Z_n^3)=9n-8$.
\end{conj*}

A significant recent progress toward this conjecture is due to Zakharov~\cite{Zakharov}, who studied the problem for prime moduli using a convex-geometric approach. Let $p$ be an odd prime. A set $A\subseteq\Z_p^r$ is called \emph{$p$-hollow} if every relation $\sum_{v\in A}a_vv=0$, with $a_v\in\Z_{\ge0}$ and $\sum_{v\in A}a_v=p$, has $a_v=p$ for some $v\in A$. Let $w(\Z_p^r)$ denote the maximum cardinality of a $p$-hollow set in $\Z_p^r$. Taking $p-1$ copies of each element of a $p$-hollow set of maximum size gives $\mathsf s(\Z_p^r)\ge w(\Z_p^r)(p-1)+1$. Zakharov~\cite{Zakharov} proved that, for every fixed $r$, $\mathsf s(\Z_p^r)=w(\Z_p^r)p+o_r(p)$ and $w(\Z_p^r)\le\binom{2r-1}{r}+1$. For $r=3$, Zakharov \cite{Zakharov} further proved that $w(\Z_p^3)=9$ for all sufficiently large primes $p$. Consequently, $\mathsf s(\Z_p^3)=9p+o(p)$ as $p\to\infty$ through the primes. Therefore, the conjectured leading term is known asymptotically for prime moduli. However, extending this result to arbitrary odd moduli is not immediate. 
Our first result establishes this extension. 

\begin{theorem}\label{thm:main}
We have $\mathsf s(\Z_n^3)\le9n+o(n)$. In particular, for odd $n$,
\[
\mathsf s(\Z_n^3)=9n+o(n).
\]
\end{theorem}

More generally, together with Zakharov's bound on $w(\Z_p^r)$, we also prove an explicit asymptotic upper bound for $\mathsf s(\Z_n^r)$ in every fixed rank over arbitrary moduli.

\begin{theorem}\label{thm:transfer}
Let $r,W\ge1$ be fixed integers. If $w(\Z_p^r)\le W$ for all sufficiently large primes $p$, then $\mathsf s(\Z_n^r)\le Wn+o_{r,W}(n)$. In particular, for every fixed $r\ge1$, we have
\[
\mathsf s(\Z_n^r)\le
\left(\binom{2r-1}{r}+1\right)n+o_r(n).
\]
\end{theorem}

Our approach combines Alon--Dubiner-type spectral expansion~\cite{AD}, Zakharov's convex-geometric framework~\cite{Zakharov}, and restricted-sum ideas developed by Han and Ren~\cite{HanRen}. We refer to~\cite{GG} for a survey on zero-sum theory and to~\cite{Bhowmick2026,EEGKR,ElsholtzEven,FGZ,FS,GHZ,GS,GZ,Naslund,SauZak,SP,WZD,WZEta} for recent studies of the Erd\H{o}s--Ginzburg--Ziv constant of finite abelian groups.

The paper is organized as follows.
In Section~\ref{sec:prelim}, we introduce the notation and provide some auxiliary results. In Section~\ref{sec:decomposition}, we decompose a sequence into finitely many families whose restricted sumsets can be prescribed at every interior length. In \ref{sec:recombination} we develop the recombination methods in general rank. Finally, in Section~\ref{sec:main-proof} we prove the the main results.

	\section{Preliminaries}\label{sec:prelim}
	
	Let $\N$ be the set of positive integers. All groups are written additively, and for every $m\ge1$ we set $\Z_m^0=\{0\}$. For a subgroup $H\le\Z_{p^a}^r$, write its invariant-factor decomposition as
	$H\cong\Z_{p^{b_1}}\oplus\cdots\oplus\Z_{p^{b_d}}$, where $1\le b_1\le\cdots\le b_d\le a$, and define $r(H)=d$; in particular, $r(\{0\})=0$.

	A sequence over $\Z_n^r$ is written as $S=x_1\cdot\ldots\cdot x_\ell$, where repetitions are allowed. We write $|S|=\ell$, $\sigma(S)=\sum_{i=1}^{\ell}x_i$, and $\supp(S)=\{x_i:1\le i\le\ell\}$. We write $T\mid S$ if $T$ is a subsequence of $S$. Pairwise disjoint subsequences always mean disjoint occurrences of terms; equivalently, their product divides $S$. In counts of terms $x\in S$, multiplicities are understood. For $0\le t\le|S|$, define
	\[
	\Sigma_t(S)=\{\sigma(T):T\mid S,\ |T|=t\},
	\qquad \ssum(S)=\bigcup_{t=0}^{|S|}\Sigma_t(S),
	\]
	where $\Sigma_0(S)=\{0\}$. Thus $\ssum(S)$ includes the empty sum.

	If $M\mid n$, multiplication by $M$ induces the canonical isomorphism
	\begin{equation}\label{eq:division}
		\iota_{n,M}:\Z_{n/M}^r\longrightarrow M\Z_n^r,
		\qquad x\longmapsto Mx.
	\end{equation}
	
	\begin{definition}\label{def:block}
		Let $S$ be a sequence over $\Z_n^r$, and let $M$ be a positive divisor of $n$. An \emph{$M$-block of $S$} is a subsequence $B\mid S$ such that
		\[
		|B|=M\qquad\text{and}\qquad \sigma(B)\in M\Z_n^r.
		\]
		Its \emph{divided sum} is the uniquely determined element
		\[
		\frac{\sigma(B)}{M}:=\iota_{n,M}^{-1}(\sigma(B))\in\Z_{n/M}^r.
		\]
	\end{definition}
A finite indexed collection $(B_\alpha)_{\alpha\in A}$ of blocks is \emph{pairwise disjoint} if $\prod_{\alpha\in A}B_\alpha\mid S$. Thus disjointness concerns occurrences, including when different blocks contain equal group elements. We use \emph{block} when its length is clear from the context.

    If $q=n/M$, selecting $t$ disjoint $M$-blocks uses exactly $Mt$ original terms, and
	\begin{equation}\label{eq:mass}
		\frac{t}{q}=\frac{Mt}{n}.
	\end{equation}
	This identity keeps track of the proportion of original terms throughout our proofs.
	
	\begin{definition}\label{def:interior-complete}
		Let $q_j\to\infty$ be positive integers, let $K_j\le\Z_{q_j}^r$, and let $Y_j$ be a sequence over $\Z_{q_j}^r$ satisfying
		$\supp(Y_j)\subseteq b_j+K_j$ and $|Y_j|/q_j\longrightarrow c>0$.
		The family $(Y_j)_j$ is \emph{complete in the interior}, relative to the moduli $q_j$ and the supporting cosets $b_j+K_j$, if for every $\varepsilon>0$ there exists $j_0(\varepsilon)$ such that, for every $j\ge j_0(\varepsilon)$,
		\begin{equation}\label{eq:complete}
			\Sigma_t(Y_j)=tb_j+K_j
			\quad\text{for every integer }t\text{ satisfying }
			\varepsilon q_j\le t\le |Y_j|-\varepsilon q_j.
		\end{equation}
	\end{definition}
	
	It is important that the same threshold $j_0(\varepsilon)$ work for every indicated $t$. The definition is independent of the representative $b_j$ of its coset and is preserved on subsequences. In particular, it applies whenever $t_j/q_j\to\theta\in(0,c)$. A family of constant sequences with these length limits is complete in the interior relative to $K_j=\{0\}$. Informally, the property allows every sum compatible with the supporting coset whenever both the selected and the unselected terms occupy a positive proportion of the modulus.

    Based on the result of Alon-Dubiner~\cite{AD} and the authors' recent result \eqref{eq:uniform-background}, there is a constant $C_r$ such that
	\begin{equation}\label{eq:linear-egz-bound}
		\mathsf s(\Z_m^d)\le C_rm
		\qquad(m\ge1,\ 0\le d\le r).
	\end{equation}
 The case $d=0$ follows from the convention $\mathsf s(\Z_m^0)=1$. 
	
	The next lemma is the standard inductive method in zero-sum theory; see, for example, \cite[Lemma~2.1(2)]{FGZ} and \cite[Equation~(2)]{AD}.
	
	\begin{lemma}\label{lem:multiplicative}
		For all $Q,m\in\N$,
		\begin{equation}\label{eq:multiplicative}
			\mathsf s(\Z_{Qm}^r)\le m\bigl(\mathsf s(\Z_Q^r)-1\bigr)+\mathsf s(\Z_m^r).
		\end{equation}
	\end{lemma}

The following result, due to Zakharov, will be used in the proof of Theorem \ref{thm:main}.

\begin{lemma}\label{prop:weak-nine}
For prime $p$ and positive integer $r$, $w(\Z_p^r)\le\binom{2r-1}{r}+1$.
For every sufficiently large prime $p$, one has $w(\Z_p^3)=9$.
\end{lemma}
\begin{proof}
    The first result is proved in \cite{Zakharov}. The second one is proved in the version 4 of \cite{Zakharov}.
\end{proof}
	
	\section{Decomposition into disjoint blocks}\label{sec:decomposition}
	
    In this section, we decompose a sequence into finitely many families whose restricted sumsets can be prescribed at every interior length. 

	\begin{proposition}\label{prop:decomposition}
		Let $N_j\to\infty$ be prime powers, and let $S_j$ be a sequence over $\Z_{N_j}^r$ satisfying $|S_j|/N_j\to c>0$. For every $\delta>0$, there is a subsequence, relabeled by $j$, and a finite nonempty index set $E$, independent of $j$, with the following properties. For every $i\in E$ and every $j$, there are a divisor $M_{i,j}\mid N_j$, an integer $\ell_{i,j}\ge1$, and $M_{i,j}$-blocks
		\[
		B_{i,j,1},\ldots,B_{i,j,\ell_{i,j}}\mid S_j.
		\]
		For each fixed $j$, the collection of all these blocks, over every $i\in E$, is pairwise disjoint. Moreover,
		\begin{gather*}
			q_{i,j}=N_j/M_{i,j}\longrightarrow\infty,\qquad
			Y_{i,j}=\prod_{\nu=1}^{\ell_{i,j}}
			\frac{\sigma(B_{i,j,\nu})}{M_{i,j}}
			,\\
			\supp(Y_{i,j})\subseteq b_{i,j}+K_{i,j},\qquad
			K_{i,j}\le\Z_{q_{i,j}}^r.
		\end{gather*}
		For every $i\in E$, the family $(Y_{i,j})_j$ is complete in the interior relative to $q_{i,j}$ and $b_{i,j}+K_{i,j}$, and
		\[
		\frac{|Y_{i,j}|}{q_{i,j}}\longrightarrow c_i>0,
		\qquad \sum_{i\in E}c_i\ge(1-\delta)c.
		\]
	\end{proposition}
	
	The product defining $Y_{i,j}$ denotes concatenation of the divided sums, so $Y_{i,j}$ is a sequence over $\Z_{q_{i,j}}^r$. The index $i$ labels a family of blocks across the chosen subsequence; the index $\nu$ labels individual blocks within that family for a fixed $j$. In particular, $|E|$ is fixed, although the numbers $\ell_{i,j}$ and the block lengths $M_{i,j}$ may vary with $j$. By \eqref{eq:mass}, $|Y_{i,j}|/q_{i,j}=M_{i,j}\ell_{i,j}/N_j$ is precisely the normalized number of original terms retained by family $i$. We shall prove Proposition \ref{prop:decomposition} at the end of this section, after the following preparations.
	
    For $N\in\N$ and a real number $K\ge0$, define the cyclic interval
	\begin{equation}\label{eq:cyclic-interval}
		I_K^{(N)}=\{a+N\Z:a\in\Z,\ |a|\le K\}\subseteq\Z_N.
	\end{equation}
	Equivalently, if $x=a+N\Z$, then $x\in I_K^{(N)}$ if and only if $\min_{z\in\Z}|a-zN|\le K$. We abbreviate this set to $I_K$ only when the modulus has been fixed; in a statement indexed by $j$, it means $I_K^{(N_j)}$. If $2K<N$, its elements have unique representatives in $[-K,K]\cap\Z$ and $|I_K^{(N)}|=2\lfloor K\rfloor+1$. A translate $a+I_K^{(N)}$ is always taken in $\Z_N$, so it may cross the residue $0$.
	
	The following proof is motivated by the strategy of \cite[Proposition~2.4 and Section~2.3]{AD}, together with a simple probabilistic argument. The next estimate applies to every nonzero homomorphism of the supporting subgroup, including those with a proper image in $\Z_N$. 
	
	\begin{lemma}\label{lem:growth}
		Let $N$ be a prime power, let $\{0\}\ne H\le \Z_N^r$, and let $\varnothing\ne Z\subseteq H$ with $s=|Z|\le h/2$, where $h=|H|$. Suppose that $A$ is a nonempty sequence over $H$ and
		\begin{equation}\label{eq:linear-thickness}
			|\{u\in A:\xi(u)\notin I_K\}|\ge\varepsilon|A|
		\end{equation}
		for every nonzero homomorphism $\xi:H\to \Z_N$, where $0<\varepsilon\le1$, $K\ge1$, and $N\ge4K$. Then there exists $u$ of $A$ such that
		\begin{equation}\label{eq:growth}
			|(Z+u)\setminus Z|\ge
			\frac1{4rN}\max\{h^{1/r}s^{1-1/r},K\varepsilon s\}.
		\end{equation}
	\end{lemma}
	
	\begin{proof}
		Define $\beta(v)=|(Z+v)\setminus Z|$ and $B=\max_{u\in\supp(A)}\beta(u)$. Translation invariance gives
		$\beta(v+w)\le\beta(v)+\beta(w)$.
		We prove two lower bounds for $B$.
		
		First, $\supp(A)$ generates $H$. Otherwise a nonzero homomorphism $H\to \Z_N$ vanishing on the subgroup generated by $A$ would contradict \eqref{eq:linear-thickness}. Write $N=p^a$. Since $r(H)\le r$, choose at most $r$ terms of $A$ whose images form a basis of $H/pH$. It is easy to see that these terms generate $H$. Repeating a chosen term if necessary, we obtain a surjection $\pi:\Z_N^r\longrightarrow H$ with $\pi(e_i)=u_i$ for $1\le i\le r$ and fibers of size $F=N^r/h$. Set $L=\left\lceil N(2s/h)^{1/r}\right\rceil$. Then $2\le L\le N$. If $t$ is uniform on $\{0,\ldots,L-1\}^r$, every value of $\pi(t)$ has probability at most $F/L^r\le1/(2s)$. Consequently,
		\[
		\mathbb E|Z\cap(Z+\pi(t))|
		=\sum_{x,y\in Z}\Pr(\pi(t)=x-y)\le\frac s2,
		\]
		and hence by the linearity of expectation, we have
		\[
		\frac s2\le\mathbb E\beta(\pi(t))
		\le\sum_{i=1}^r\mathbb E(t_i)\beta(u_i)
		\le\frac{r(L-1)}2 B.
		\]
		Since $L-1<N(2s/h)^{1/r}$, we obtain
		\begin{equation}\label{eq:growth-small}
			B\ge\frac{h^{1/r}s^{1-1/r}}{r2^{1/r}N}\ge\frac{h^{1/r}s^{1-1/r}}{2rN}.
		\end{equation}
		
		For the second bound, let $M=\lceil N/K\rceil$ and let $U=tu$, where $u$ is uniform among the occurrences of $A$ and $t$ is independently uniform on $\{1,\ldots,M\}$. Every character of $H$ has the form
		$\chi(u)=\exp(2\pi i\xi(u)/N)$
		for a homomorphism $\xi:H\to \Z_N$. For a nontrivial character, at least an $\varepsilon$ proportion of the terms $u$ have a representative $d\in(-N/2,N/2]$ of $\xi(u)$ with $|d|>K$. For these terms, the geometric-sum estimate gives
		$\left|\frac1M\sum_{t=1}^M\chi(tu)\right|\le\frac{N}{2M|d|}\le\frac12$.
		Thus $\operatorname{Re}\mathbb E\chi(U)\le1-\varepsilon/2$. Let
		$\widehat{1_Z}(\chi)=h^{-1/2}\sum_{z\in Z}\overline{\chi(z)}$ and let $\widehat H$ be the character group of $H$. By the Fourier orthogonality, we have
		$$
			\mathbb E\beta(U)
			=\sum_{\chi\in\widehat H}|\widehat{1_Z}(\chi)|^2
			\bigl(1-\operatorname{Re}\mathbb E\chi(U)\bigr)
			\ge\frac\varepsilon2\left(s-\frac{s^2}{h}\right)
			\ge\frac{\varepsilon s}{4}.$$
		On the other hand, $\mathbb E\beta(U)\le(M+1)B/2$. Since $N\ge4K$, we have $M+1\le3N/(2K)$, and therefore
		$B\ge K\varepsilon s/(3N)$.
		Together with \eqref{eq:growth-small}, this proves \eqref{eq:growth}.
	\end{proof}
	
	\begin{lemma}\label{lem:pairs}
		Let $N_j\to\infty$ be prime powers, let $H_j\le \Z_{N_j}^r$, and let $T_j$ be a sequence over $\Z_{N_j}^r$ satisfying $\supp(T_j)\subseteq b_j+H_j$ and $|T_j|\ge\mu N_j$ for some $\mu>0$. Suppose that there is a fixed $0<\varepsilon\le1$ such that, for every fixed $K\in\N$, for all sufficiently large $j$,
		\begin{equation}\label{eq:affine-thickness}
			|\{x\in T_j:\xi(x-b_j)\notin a+I_K\}|\ge\varepsilon|T_j|
		\end{equation}
		for every nonzero homomorphism $\xi:H_j\to \Z_{N_j}$ and every $a\in \Z_{N_j}$. Then, for all sufficiently large $j$, there are $R_j=o(N_j)$ pairwise disjoint pairs of terms $(a_1,b_1),\ldots,(a_{R_j},b_{R_j})$ of $T_j$ such that
		\[
		\ssum\bigl((b_1-a_1)\cdot\ldots\cdot(b_{R_j}-a_{R_j})\bigr)=H_j.
		\]
	\end{lemma}
	
	\begin{proof}
		If $H_j=\{0\}$, no pairs are needed. We may therefore assume that $H_j\ne\{0\}$, and, for simplicity, we omit the subscript \(j\) in the following argument. Fix a large $K$.
        The condition \eqref{eq:affine-thickness} is preserved, with \(\varepsilon/2\) in place of \(\varepsilon\), as long as at most \(\varepsilon\mu N/4\) occurrences have been removed. Indeed, if \(X\) is the remaining sequence, then \(|T|-|X|\le\varepsilon\mu N/4\le\varepsilon|T|/4\). Hence, for every nonzero \(\xi:H\to\mathbb Z_N\) and every \(a\in\mathbb Z_N\),
$$
|\{x\in X:\xi(x-b)\notin a+I_K\}|
\ge\varepsilon|T|-(|T|-|X|)
\ge\frac{3\varepsilon}{4}|T|
\ge\frac{\varepsilon}{2}|X|.
$$
This allows us to apply Lemma \ref{lem:growth} repeatedly while constructing the disjoint pairs.
		
		Start with $Z=\{0\}$. Fix any remaining occurrence $a\in X$ and apply Lemma~\ref{lem:growth} to the sequence of differences $b-a$, where $b$ runs through $X$. Condition \eqref{eq:linear-thickness} follows from \eqref{eq:affine-thickness}, by shifting the center of the strip. As long as $s=|Z|\le h/2$, where $h=|H|$, we can choose $b$ so that, with $u=b-a$ and $Z'=Z+\{0,u\}$,
		\begin{equation}\label{eq:pair-growth}
			|Z'|-|Z|\ge\frac1{8rN}
			\max\{h^{1/r}s^{1-1/r},K\varepsilon s\}.
		\end{equation}
		The difference $u$ is nonzero. Hence $a$ and $b$ are distinct occurrences, and we delete both before choosing the next pair.
		
		We estimate the number of steps. Define $A_0=K\varepsilon$ and take $K$ large enough that $A_0\ge2$. Set $x=N(s/h)^{1/r}\ge1$. The first bound in \eqref{eq:pair-growth} implies
		\[
		x'\ge\left(x^r+\frac{x^{r-1}}{8r}\right)^{1/r}
		\ge x+\kappa_r
		\]
		for a constant $\kappa_r>0$ depending only on $r$. For example, one may take $\kappa_r=1/(16r^2)$, by the binomial theorem and $x\ge1$.
		Thus after $O_r(N/A_0+1)$ steps, either $s>h/2$ or $s\ge h/A_0^r$. Thereafter the second bound gives
		$s'\ge\left(1+\frac{A_0}{8rN}\right)s$.
		For fixed $K$ and sufficiently large $N$, a further $O_r((N/A_0)\log A_0+1)$ steps suffice to obtain $s>h/2$.
		
		Repeat the construction with fresh endpoints. This gives two subset-sum sets of size greater than $h/2$. Clearly, their sum is $H$. The total number of pairs is		$O_{r,\varepsilon}\!\left(\frac{N\log K}{K}+1\right)$.
		Choose $K$ sufficiently large and then $N$ sufficiently large so that this number is less than $\varepsilon\mu N/8$. At most twice as many occurrences are deleted, so the assumed deletion bound holds throughout both constructions.
		
		For each sufficiently large $j$, choose a family with the minimum possible number $R_j$ of pairs. For every sufficiently large fixed $K$, the preceding construction gives
		\[
		\limsup_{j\to\infty}\frac{R_j}{N_j}
		\le C_{r,\varepsilon}\frac{\log K}{K}.
		\]
		As $K\to\infty$, we have $R_j=o(N_j)$.
	\end{proof}
	
	\begin{corollary}\label{cor:complete}
		Under the hypotheses of Lemma~\ref{lem:pairs}, suppose also that $|T_j|/N_j\to c>0$. Then $(T_j)_j$ is complete in the interior relative to the moduli $N_j$ and the supporting cosets $b_j+H_j$.
	\end{corollary}
	
	\begin{proof}
By Lemma~\ref{lem:pairs}, there are $R_j=o(N_j)$ pairwise
disjoint pairs
\(
(u_1,v_1),\ldots,(u_{R_j},v_{R_j})
\)
of terms of $T_j$ such that
\(
\Sigma^0\bigl((v_1-u_1)\cdots(v_{R_j}-u_{R_j})\bigr)
=H_j\).
Choosing exactly one endpoint from each pair gives
$R_j$ terms whose sums cover the coset
$R_jb_j+H_j$. Now let $t$ be an integer satisfying
\(
R_j\le t\le |T_j|-R_j
\).
Since the pairs use $2R_j$ occurrences, at least
$t-R_j$ occurrences remain outside them. Choose
and fix any $t-R_j$ such occurrences, and let
$s$ denote their sum. Since
$
s\in(t-R_j)b_j+H_j$, adding one endpoint from each pair gives
$t$-term subsequences whose sums cover
\(
s+(R_jb_j+H_j)=tb_j+H_j\).
The reverse inclusion follows from
$\operatorname{supp}(T_j)\subseteq b_j+H_j$.
Consequently, we have
$
\Sigma_t(T_j)=tb_j+H_j$, where
\(R_j\le t\le |T_j|-R_j\). Finally, since $R_j=o(N_j)$, for every fixed
$\varepsilon>0$ and all sufficiently large $j$,
we have $R_j\le\varepsilon N_j$. Thus the
preceding equality holds for every integer
$
\varepsilon N_j\le t\le
|T_j|-\varepsilon N_j,
$
which is precisely the completeness property
in \eqref{eq:complete}.
	\end{proof}
	
	We next handle concentration in a subgroup. The relevant measure of progress is its rank after multiplication by a suitable divisor, rather than its index alone.
	
	\begin{lemma}\label{lem:quotient}
		Let $N$ be a prime power, let $H\le \Z_N^r$, let $\supp(S)\subseteq b+H$, and let $m\mid N$. All but fewer than $C_rm$ terms of $S$ can be partitioned into blocks of length $m$ whose divided sums lie in $\bar b+K$, where
		\[
		\bar b\in \Z_{N/m}^r\text{ is the reduction of }b,
		\qquad K\le \Z_{N/m}^r,\qquad K\cong mH.
		\]
		If $H\cong \Z_{d_1}\oplus\cdots\oplus \Z_{d_d}$, with $1<d_1\mid\cdots\mid d_d$, then
		\[
		r(K)=\#\{i:d_i>m\}.
		\]
		In particular, if $d_1=o(N)$ along a sequence of such groups, choosing $m=d_1$ lowers the rank, gives $N/m\to\infty$, and preserves the limiting ratio of sequence length to modulus.
	\end{lemma}
	
	\begin{proof}
		The quotient $H/mH$ has exponent dividing $m$ and at most $r$ generators, so choose an injective homomorphism $H/mH\hookrightarrow\Z_m^r$. View the images of the deviations $x-b$ as a sequence in this ambient group. Applying \eqref{eq:linear-egz-bound} to $\Z_m^r$ gives a zero-sum subsequence of exactly $m$ terms whenever at least $C_rm$ terms remain. Injectivity makes its sum zero in $H/mH$ as well. Repeating the selection leaves fewer than $C_rm$ terms. This use of the ambient group is necessary because the exponent of $H/mH$ may be smaller than $m$. If $B$ is one of the corresponding blocks, then
		$\sigma(B)\in mb+mH\subseteq m\Z_N^r$.
		Under \eqref{eq:division}, its divided sum belongs to $\bar b+K$, where $K$ is the inverse image of $mH$. In particular, $K\cong mH$.
		
		Since all relevant orders are powers of the same prime,
		$mH\cong\bigoplus_{i=1}^d \Z_{d_i/\gcd(d_i,m)}$.
		A factor is nontrivial exactly when $d_i>m$, proving the rank formula. If $R<C_rm$ terms remain, $q=N/m$, and $Y$ is the sequence of divided sums, then
		\begin{equation}\label{eq:quotient-mass}
			\frac{|Y|}{q}=\frac{|S|-R}{N}.
		\end{equation}
		Thus $m=o(N)$ makes the loss in the normalized length tend to zero.
	\end{proof}
	
	The reduction in Lemma~\ref{lem:quotient} can be iterated. Suppose that each current term $x\in \Z_N^r$ represents a block of $P$ original terms in $\Z_n^r$, where $n=PN$, and that the sum of this original block is $Px$. If $m$ current terms have sum $my$ under \eqref{eq:division}, the product of their original blocks has length $Pm$ and sum $Pmy$. Consequently, disjointness, divisibility of block lengths, and \eqref{eq:mass} are preserved at every stage.
	
	\begin{proof}[Proof of Proposition~\ref{prop:decomposition}]
		We prove the stronger assertion for sequences supported in cosets of subgroups of rank at most $d$, by induction on $0\le d\le r$. The observation following Lemma~\ref{lem:quotient} allows current terms to represent previously selected blocks. It is enough to consider $0<\delta<1$.
		
		For $d=0$, the sequence is constant. Take each occurrence as a $1$-block, use one family, and choose the supporting subgroup $\{0\}$; this family is complete in the interior. Suppose that the assertion is known for smaller ranks. Passing to a subsequence, we may assume that $\supp(S_j)\subseteq a_j+H_j$ and $r(H_j)=d$. Choose $0<\eta<1$ with $(1-\eta)^2\ge1-\delta$, and set $\alpha=(1-\eta)c$.
		
		For each sufficiently large $j$, choose a coset $C_j=b_j+L_j\subseteq a_j+H_j$ maximizing
		\begin{equation}\label{eq:excess}
			\frac{|(S_j)_{x+L}|-\alpha N_j}{|L|}
			\qquad(L\le H_j,\ x\in a_j+H_j).
		\end{equation}
		Let $Y_j=(S_j)_{C_j}$, $D_j=[H_j:L_j]$, and $c'_j=|Y_j|/N_j$. Comparison with the whole coset $a_j+H_j$ gives
		\begin{equation}\label{eq:excess-comparison}
			D_j(c'_j-\alpha)\ge |S_j|/N_j-\alpha>0.
		\end{equation}
		Pass to a subsequence on which $c'_j\to c'\ge\alpha$, the rank of $L_j$ is fixed, and $D_j$ is either constant or tends to infinity.
		
		If $r(L_j)<d$, apply induction to $Y_j$ with tolerance $\eta$. The retained normalized length is at least $(1-\eta)c'\ge(1-\eta)^2c\ge(1-\delta)c$.
		If $r(L_j)=d$ and $D_j\to\infty$, let $d_j$ be the smallest invariant factor of $L_j$. Then
		\begin{equation}\label{eq:small-factor}
			d_j^d\le |L_j|=\frac{|H_j|}{D_j}\le\frac{N_j^d}{D_j},
			\qquad \frac{d_j}{N_j}\longrightarrow0.
		\end{equation}
		Lemma~\ref{lem:quotient} lowers the rank with a negligible loss, and induction with tolerance $\eta$ again applies.
		
		It remains to consider $r(L_j)=d$ and $D_j=D$ constant. By \eqref{eq:excess-comparison}, $D(c'-\alpha)\ge c-\alpha>0$. If $C'\subsetneq C_j$ is a coset of a proper subgroup $L'<L_j$, maximality in \eqref{eq:excess} gives
		\begin{equation}\label{eq:proper-subcoset}
			\frac{|(Y_j)_{C'}|}{N_j}
			\le\alpha+\frac{c'_j-\alpha}{[L_j:L']}
			\le\frac{c'_j+\alpha}{2}.
		\end{equation}
		Choose $0<\varepsilon<1$ such that $(1-\varepsilon)c'>(c'+\alpha)/2$.
		The same strict inequality holds with $c'_j$ in place of $c'$ for all sufficiently large $j$.
		
		If $Y_j$ satisfies \eqref{eq:affine-thickness} for every fixed $K$, Corollary~\ref{cor:complete} provides a terminal family that is complete in the interior. Otherwise, for some fixed $K$, along a subsequence there are nonzero homomorphisms $\xi_j:L_j\to \Z_{N_j}$ and centers $u_j\in \Z_{N_j}$ for which
		\[
		|\{x\in Y_j:\xi_j(x-b_j)\in u_j+I_K\}|
		>(1-\varepsilon)|Y_j|.
		\]
		Pass to a further subsequence on which $t_j=|\im\xi_j|$ is bounded or tends to infinity. The bounded case is impossible. Indeed, distinct points of $\im\xi_j$ are separated by $N_j/t_j>2K$, so the strip meets at most one fiber of $\xi_j$. This is a coset of the proper subgroup $\ker\xi_j$, contradicting \eqref{eq:proper-subcoset}.
		
		Hence $t_j\to\infty$. Split the strip into its at most $2K+1$ kernel fibers. Label these fibers by the offsets $-K,\ldots,K$, allowing empty fibers. After passing to a subsequence, each of their normalized lengths converges and $r(\ker\xi_j)$ is constant. Discard fibers whose limiting normalized length is zero. The sum of the positive limits is at least $(1-\varepsilon)c'>\alpha$. The subgroup $\ker\xi_j$ either has rank less than $d$, or has rank $d$ and index at least $t_j$ in $H_j$. In the latter case, \eqref{eq:small-factor} and Lemma~\ref{lem:quotient} lower the rank. Apply induction, with tolerance $\eta$, to each remaining fiber. The total retained normalized length is at least $(1-\eta)\alpha=(1-\eta)^2c\ge(1-\delta)c$.
		
		At each recursive step there are finitely many children and the supporting rank strictly decreases. Thus the recursion has a finite tree of terminal families. Label its leaves by a finite index set $E$. Because only finitely many subsequence extractions have been made, they have a common infinite subsequence, which we relabel by $j$. The set $E$ is therefore independent of $j$ and is nonempty because the retained normalized length is positive. Lemma~\ref{lem:quotient} and \eqref{eq:mass} preserve disjoint original occurrences, divisibility of block lengths, and normalized lengths. Every terminal modulus tends to infinity. This completes the proof.
	\end{proof}
	
	\section{Recombination with exact lengths}\label{sec:recombination}
	
	The block families may have different lengths, so their available proportions cannot be combined freely. We retain the exact divisibility and zero-sum conditions in finite torus subgroups, then use a minimum-volume simplex to find an admissible choice strictly inside the available capacities. This is the only step that uses an upper bound for the weak EGZ constant.
	
	\begin{proposition}\label{prop:recombination}
		Let $W\ge1$ be an integer such that $w(\Z_p^r)\le W$ for all sufficiently large primes $p$. Let $n_j\to\infty$, and suppose that sequences over $\Z_{n_j}^r$ contain finitely many families of pairwise disjoint blocks as in Proposition~\ref{prop:decomposition}, with lengths $M_{i,j}\mid n_j$, divided-sum sequences $Y_{i,j}$, and positive limits $c_i$. Assume that every family is complete in the interior and that
		\begin{equation}\label{eq:common-modulus}
			Q_j=\frac{n_j}{\lcm_{i\in E}M_{i,j}}\longrightarrow\infty.
		\end{equation}
		If $\sum_{i\in E}c_i>W$, then, for all sufficiently large $j$, the original sequence has a zero-sum subsequence of length $n_j$.
	\end{proposition}
	
	The coordinate $\theta_i$ in the argument below will represent the proportion of the $n_j$ original terms selected from family $i$. We seek proportions summing to one and lying strictly below the available proportions $c_i$. Completeness in the interior will then allow us to prescribe a sum in each supporting coset. The finite groups introduced below record both the exact length congruence and the sum congruence. Their limits allow us to choose the proportions geometrically.
	
	We shall prove Proposition~\ref{prop:recombination} by contradiction. Suppose that its conclusion fails for infinitely many indices \(j\). Thus, there exists an infinite set \(J\subseteq\mathbb N\) such that, for every \(j\in J\), the original sequence contains no zero-sum subsequence of length \(n_j\). In what follows, we restrict \(j\) to \(J\), and all limits are understood along this infinite set.

     We first observe that $c_i\le1$ for all $i$. Indeed, if $c_i>1$, property \eqref{eq:complete} with $t=q_{i,j}=n_j/M_{i,j}$ gives $\Sigma_{q_{i,j}}(Y_{i,j})=K_{i,j}$, which contains zero. The corresponding $q_{i,j}$ blocks already give $n_j$ original terms with sum zero. Set $A=\sum_{i\in E}c_i>W$.
	
	Choose a lift of $b_{i,j}\in \Z_{q_{i,j}}^r$ to $\Z_{n_j}^r$ and denote it by the same symbol. Let $H_{i,j}=M_{i,j}K_{i,j}\le\Z_{n_j}^r$, using the embedding \eqref{eq:division}. Selecting $t_i$ blocks from family $i$, or $k_i=M_{i,j}t_i$ original terms, realizes all elements of
	\begin{equation}\label{eq:sum-coset}
		k_ib_{i,j}+H_{i,j}
	\end{equation}
	whenever $t_i/q_{i,j}$ stays in a compact subinterval of $(0,c_i)$.
	
	For $\varnothing\ne I\subseteq E$, write $H_j(I)=\sum_{i\in I}H_{i,j}$, $D_j(I)=\prod_{i\in I}M_{i,j}\Z_{n_j}$, and $\T^I=(\R/\Z)^I$.
	Define a finite subgroup of $\T^I$ by
	\begin{equation}\label{eq:Gamma}
		\Gamma_j(I)=\left\{\frac{k}{n_j}\bmod\Z^I:
		\begin{array}{l}
			k\in D_j(I),\quad \sum_{i\in I}k_i=0\text{ in }\Z_{n_j},\\
			\sum_{i\in I}k_ib_{i,j}\in H_j(I)
		\end{array}\right\}.
	\end{equation}
	Here a residue $k_i\in \Z_{n_j}$ determines $k_i/n_j\in\R/\Z$. The expression is independent of the chosen lifts of $b_{i,j}$, because $k_i$ is a multiple of $M_{i,j}$ and $M_{i,j}q_{i,j}=n_j$.
	
	The nonempty closed subsets of a compact metric space form a compact space in the Hausdorff metric. Passing to a subsequence, we may therefore assume that all the finitely many groups $\Gamma_j(I)$ converge in Hausdorff distance to closed subgroups $\Gamma(I)\le\T^I$. Set
	\begin{equation}\label{eq:Lambda}
		\Lambda(I)=\{\lambda\in\Z^I:
		\lambda\cdot x=0\text{ in }\R/\Z
		\text{ for every }x\in\Gamma(I)\}.
	\end{equation}
Since $\Lambda(I)\le\mathbb Z^I$, it is a finitely
generated free abelian group. We retain the full
lattice $\Lambda(I)$, rather than merely its rational
span, to preserve the congruence information.
Recall that every continuous character of
$\mathbb T^I$ has the form
\[
x\longmapsto e^{2\pi i\lambda\cdot x},
\qquad \lambda\in\mathbb Z^I.
\]
By the double-annihilator theorem for closed
subgroups of compact abelian groups
(see \cite[Chapter~7, Proposition~38]{Morris1977}),
a closed subgroup is the common kernel of all
characters vanishing on it. Therefore,
\begin{equation}\label{eq:annihilator}
\Gamma(I)=
\left\{x\in\mathbb T^I:
\lambda\cdot x=0\text{ in }\mathbb R/\mathbb Z
\text{ for every }\lambda\in\Lambda(I)\right\}.
\end{equation}
	
	\begin{lemma}\label{lem:relations}
		For a fixed $\lambda\in\Z^I$, one has $\lambda\in\Lambda(I)$ if and only if $\lambda$ annihilates $\Gamma_j(I)$ for every sufficiently large $j$. Moreover,
		\begin{equation}\label{eq:restriction}
			\mathbf1_I\in\Lambda(I),\qquad
			\Lambda(J)|_I\subseteq\Lambda(I)
			\quad(\varnothing\ne I\subseteq J\subseteq E),
		\end{equation}
		where $\mathbf1_I$ is the vector with all coordinates equal to $1$ and $\Lambda(J)|_I$ denotes restriction of coordinates. For every $\lambda\in\Lambda(I)$ and every sufficiently large $j$, there are $t_0\in \Z_{n_j}$ and $t\in \Z_{n_j}^r$ such that
		\begin{equation}\label{eq:affine-relations}
			\lambda_i\equiv t_0+t\cdot b_{i,j}\pmod{q_{i,j}}
			\qquad(i\in I).
		\end{equation}
	\end{lemma}
	
	\begin{proof}
		If $\lambda$ annihilates $\Gamma_j(I)$ eventually, it annihilates the Hausdorff limit. Conversely, if $\lambda\in\Lambda(I)$, continuity and Hausdorff convergence show that the maximum distance of $\lambda\cdot x$ from $0$ in $\R/\Z$, for $x\in\Gamma_j(I)$, tends to zero. A nontrivial finite subgroup of $\R/\Z$ has a point at distance at least $1/3$ from zero. Hence the image of $\Gamma_j(I)$ under this character must eventually be trivial. The first assertion follows.
		
		The length condition in \eqref{eq:Gamma} gives $\mathbf1_I\in\Lambda(I)$. Extending coordinates by zero embeds $\Gamma_j(I)$ into $\Gamma_j(J)$, so stabilization gives the restriction assertion in \eqref{eq:restriction}.
		
		For the last assertion, suppress $j$ and consider the homomorphism $F:D(I)\longrightarrow \Z_n\oplus(\Z_n^r/H(I))$ given by $F(k)=\left(\sum_{i\in I}k_i,\sum_{i\in I}k_ib_i+H(I)\right)$. The character $k\mapsto\sum_{i\in I}\lambda_i k_i/n\pmod\Z$ is trivial on $\ker F$. It therefore descends to $\im F$ and extends to the displayed finite abelian group. Such an extension exists because a character into $\R/\Z$ can be extended one element at a time by taking a root of its prescribed multiple. The extended character has the form
		\[
		(u,z+H(I))\longmapsto\frac{t_0u+t\cdot z}{n}\pmod\Z,
		\]
		where $t$ annihilates $H(I)$. Evaluate at a coordinate vector of value $M_i$ in $D(I)$. We obtain
		\[
		M_i\lambda_i\equiv M_i(t_0+t\cdot b_i)\pmod n,
		\]
		which is precisely \eqref{eq:affine-relations}.
	\end{proof}
	
	The next lemma transfers the finite collection of congruence relations to an arbitrarily large prime field. 

\begin{lemma}\label{lem:specialization}
Let $f_1,\ldots,f_m\in\mathbb Z[X_1,\ldots,X_d]$.
If the system
\[
f_1=\cdots=f_m=0
\]
has a solution modulo arbitrarily large positive integers,
then it has a solution over $\mathbb Z_p$ for infinitely
many primes $p$.
\end{lemma}

\begin{proof}
Let
\[
B=\mathbb Z[X_1,\ldots,X_d]/(f_1,\ldots,f_m).
\]
The hypothesis implies that no nonzero integer vanishes
in $B$: otherwise, evaluating at solutions modulo
unbounded integers $N$ would force $N$ to divide
a fixed nonzero integer. Hence
\[
B\otimes_{\mathbb Z}\mathbb Q\ne0.
\]
By Zariski's lemma, this algebra has a maximal-ideal
quotient which is a finite extension of $\mathbb Q$.
Thus the system has a common zero over
$\overline{\mathbb Q}$.

It is known that any finite system of polynomial
equations over $\mathbb Z$ with a common zero over
$\overline{\mathbb Q}$ is solvable over $\mathbb Z_p$
for infinitely many primes $p$; see
\cite[Proposition~2.7]{vandenDries1991} and
\cite[Theorem~1.2 and the discussion following it]
{Jarviniemi2022}.
This proves the lemma.
\end{proof}

	Let $\Delta_E=\{x\in\R_{\ge0}^E:\sum_{i\in E}x_i=1\}$ and, for $x\in\Delta_E$, let $\supp(x)=\{i\in E:x_i>0\}$.
	Define
	\begin{equation}\label{eq:admissible}
		\cS=\{\theta\in\Delta_E:
		\lambda\cdot\theta_I\in\Z
		\text{ for every }\lambda\in\Lambda(I),
		\ I=\supp(\theta)\}.
	\end{equation}
	Here $\theta_I$ is the restriction to the indicated coordinates. Every coordinate vertex of $\Delta_E$ belongs to $\cS$.
	
	\begin{lemma}\label{lem:capacities}
		There is a vector $\theta\in\cS$ such that $\theta_i<c_i$ for every $i\in E$.
	\end{lemma}
	
	\begin{proof}
		Define $\mathcal T=\overline{\cS}$ and $b=(c_i/A)_{i\in E}$. Choose a simplex $P=\conv\{x_1,\ldots,x_m\}$ with $x_\nu\in\mathcal T$, containing $b$ and first minimizing $m$. Such a simplex exists because $\cS$ contains the coordinate vertices. If $m=1$, then $b\in\mathcal T$ and $b_i=c_i/A<c_i$ for all $i$, since $A>W\ge1$. A sufficiently close point of $\cS$ proves the assertion.
		
		Suppose that $m\ge2$. Among the feasible $m$-tuples, minimize the $(m-1)$-dimensional Euclidean volume. This minimum exists: the set of pairs consisting of an $m$-tuple in $\mathcal T^m$ and barycentric coefficients representing $b$ is compact, and so is its projection onto $\mathcal T^m$. Minimality of $m$ excludes affine dependence and forces every barycentric coefficient to be positive. Volume is continuous on this compact feasible set. Write
		\begin{equation}\label{eq:barycentric}
			b=\sum_{\nu=1}^m\gamma_\nu x_\nu,
			\qquad \gamma_\nu>0,\qquad\sum_{\nu=1}^m\gamma_\nu=1.
		\end{equation}
		
		We next choose rational representatives while retaining the supports of approximating points. For each $x_\nu$, take a sequence in $\cS$ converging to $x_\nu$. Passing to a subsequence fixes its support $I_\nu$ and the values of the pairings with a $\Z$-basis of $\Lambda(I_\nu)$: these values are bounded integers. The resulting affine equations have rational coefficients and a solution positive on $I_\nu$. Rational points are dense in their affine solution space. We may therefore choose $z_\nu\in\cS\cap\Q^E$ satisfying
		\begin{equation}\label{eq:rational-representatives}
			\begin{gathered}
				\supp(x_\nu)\subseteq I_\nu=\supp(z_\nu),\\
				\lambda\cdot(x_\nu)_{I_\nu}
				=\lambda\cdot(z_\nu)_{I_\nu}\in\Z
				\qquad(\lambda\in\Lambda(I_\nu)).
			\end{gathered}
		\end{equation}
		The equation $\sum_i(z_\nu)_i=1$ is included because $\mathbf1_{I_\nu}\in\Lambda(I_\nu)$. These rational points are used only for arithmetic specialization; the minimizing simplex and its barycentric identity still use the original vertices $x_\nu$.
		
		For each nonempty $I\subseteq E$, choose a $\Z$-basis of $\Lambda(I)$. For every basis vector $\lambda$, introduce an affine function $f_{I,\lambda}$ in $r$ variables, and introduce vectors $w_i$ with $r$ coordinates, for $i\in E$. Impose the equations
		\begin{equation}\label{eq:polynomial-system}
			f_{I,\lambda}(w_i)=\lambda_i\qquad(i\in I).
		\end{equation}
		The unknowns are the coordinates of the $w_i$ and the coefficients of the affine functions, so this is a finite system of integer polynomial equations. There are finitely many nonempty subsets $I\subseteq E$ and finitely many chosen basis vectors for each $\Lambda(I)$; hence the eventual assertions of Lemma~\ref{lem:relations} hold simultaneously. The system therefore has solutions modulo $Q_j$ for all sufficiently large $j$: indeed, $Q_j$ divides every $q_{i,j}$, and the lifts $b_{i,j}$ give a simultaneous choice for $w_i$. Lemma~\ref{lem:specialization} supplies a solution over $\Z_p$, where $p$ is sufficiently large that $w(\Z_p^r)\le W$ and avoids all denominators of the coordinates of the $z_\nu$. Integer linear combinations define $f_{I,\lambda}$ for all $\lambda\in\Lambda(I)$. Set $u_\nu=\sum_{i\in E}(z_\nu)_i w_i\in\Z_p^r$ for $1\le\nu\le m$, where rational coefficients are reduced modulo $p$.
		
		We claim that the $u_\nu$ are distinct and form a $p$-hollow set. Otherwise there are integers $0\le a_\nu<p$ with
		\begin{equation}\label{eq:forbidden}
			\sum_{\nu=1}^m a_\nu=p,
			\qquad \sum_{\nu=1}^m a_\nu u_\nu=0.
		\end{equation}
		Coincident points are included: coefficients $1$ and $p-1$ at two such points give this relation. Set $y=p^{-1}\sum_\nu a_\nu x_\nu$, $z=p^{-1}\sum_\nu a_\nu z_\nu$, and $J=\bigcup_{a_\nu>0}I_\nu$.
		For $\lambda\in\Lambda(J)$ and every active $\nu$, the restriction property \eqref{eq:restriction}, \eqref{eq:rational-representatives}, and the affine nature of $f_{J,\lambda}$ give
		\begin{equation}\label{eq:affine-values}
			f_{J,\lambda}(u_\nu)
			=\overline{\lambda\cdot(z_\nu)_J}
			=\overline{\lambda\cdot(x_\nu)_J},
		\end{equation}
		where the bars denote reduction of integer values modulo $p$. Applying the affine function to \eqref{eq:forbidden} shows that
		\[
		\sum_\nu a_\nu\lambda\cdot(z_\nu)_J\equiv0\pmod p.
		\]
		The constant term of the affine function disappears because $\sum_\nu a_\nu=p$. It follows that
		\[
		\lambda\cdot y_J=\lambda\cdot z_J\in\Z
		\qquad(\lambda\in\Lambda(J)).
		\]
		Now $\supp(z)=J$ and $\supp(y)\subseteq J$. Hence $(1-t)y+tz\in\cS$ for every $0<t\le1$, and consequently $y\in\mathcal T$. This interpolation is needed because $y$ may lie on a smaller coordinate face.
		
		Write $\rho_\nu=a_\nu/p$, and choose $h$ minimizing $\gamma_\nu/\rho_\nu$ among the indices with $\rho_\nu>0$. Let $t=\gamma_h/\rho_h$. Then
		\[
		b=t y+\sum_{\nu\ne h}(\gamma_\nu-t\rho_\nu)x_\nu,
		\]
		with nonnegative coefficients summing to $1$. Thus the simplex obtained by replacing $x_h$ with $y$ still contains $b$. The opposite face is unchanged, and the altitude of $y=\sum_\nu\rho_\nu x_\nu$ above that face is $\rho_h$ times the altitude of $x_h$. Its volume is therefore $\rho_h\vol(P)$, which is positive and smaller than $\vol(P)$, because $0<\rho_h<1$. This contradicts the choice of $P$ and proves the claim.
		
		Therefore $m\le w(\Z_p^r)\le W$. Some $\gamma_h\ge1/m\ge1/W$. The barycentric identity expresses the normalized capacity vector $b=(c_i/A)_{i\in E}$ as a weighted average of at most $W$ admissible limit points. A vertex with weight at least $1/W$ must therefore fit strictly below every capacity. More precisely, by nonnegativity and \eqref{eq:barycentric},
		\[
		(x_h)_i\le\frac{b_i}{\gamma_h}
		\le\frac{Wc_i}{A}<c_i\qquad(i\in E).
		\]
		A sufficiently close point of $\cS$ satisfies the same strict inequalities. This completes the proof.
	\end{proof}
	
	\begin{proof}[Proof of Proposition~\ref{prop:recombination}]
		Recall that we have restricted to an infinite set
$J\subseteq\mathbb N$ such that, for every $j\in J$,
the original sequence contains no zero-sum
subsequence of length $n_j$.
We now derive a contradiction. Choose $\theta\in\cS$ as in Lemma~\ref{lem:capacities}, and let $I=\supp(\theta)$. Then
		\[
		0<\theta_i<c_i\le1\quad(i\in I),\qquad
		\sum_{i\in I}\theta_i=1.
		\]
		By \eqref{eq:admissible} and \eqref{eq:annihilator}, $\theta_I\bmod\Z^I\in\Gamma(I)$. Hausdorff convergence gives elements of $\Gamma_j(I)$ approaching this point. Write them as $k_j/n_j\bmod\Z^I$, with standard representatives $0\le k_{i,j}<n_j$. Each $k_{i,j}$ is divisible by $M_{i,j}$. Since $0<\theta_i<1$, there is no ambiguity at the endpoints of the unit interval, and $k_{i,j}/n_j\longrightarrow\theta_i$ for every $i\in I$.
		By the length congruence in \eqref{eq:Gamma}, $\sum_{i\in I}k_{i,j}/n_j$ is an integer. It tends to $1$, so
		\begin{equation}\label{eq:exact-length}
			\sum_{i\in I}k_{i,j}=n_j
		\end{equation}
		for all sufficiently large $j$.
		
		Let $t_{i,j}=k_{i,j}/M_{i,j}$. By \eqref{eq:mass}, $t_{i,j}/q_{i,j}\longrightarrow\theta_i\in(0,c_i)$.
		Since $I$ is finite, choose $\varepsilon_0>0$ such that $3\varepsilon_0<\min_{i\in I}\{\theta_i,c_i-\theta_i\}$.
		The convergences above imply, for all sufficiently large $j$,
		\[
		\varepsilon_0q_{i,j}\le t_{i,j}
		\le |Y_{i,j}|-\varepsilon_0q_{i,j}
		\qquad(i\in I).
		\]
		Completeness in the interior therefore applies simultaneously to every active family.
		
		Membership in $\Gamma_j(I)$ gives
		\[
		-\sum_{i\in I}k_{i,j}b_{i,j}\in H_j(I)
		=\sum_{i\in I}H_{i,j}.
		\]
		Choose $h_{i,j}\in H_{i,j}$ whose sum equals the expression on the left. For each $i$, completeness realizes $k_{i,j}b_{i,j}+h_{i,j}$ using exactly $t_{i,j}$ blocks. These selections use disjoint original occurrences and have sum zero, contradicting our assumption since their total length is
		\[
		\sum_{i\in I}M_{i,j}t_{i,j}
		=\sum_{i\in I}k_{i,j}=n_j.\qedhere
		\]
	\end{proof}
	
	\section{Proof of the main theorems}\label{sec:main-proof}

We first combine Propositions~\ref{prop:decomposition}
and~\ref{prop:recombination} to establish the desired
asymptotic bound uniformly over prime powers.
The key observation is that the divisors of a prime
power are totally ordered by divisibility.
Consequently, the common residual modulus is the
minimum of the finitely many residual moduli and
therefore tends to infinity, as required by
Proposition~\ref{prop:recombination}.
We then extend the bound to arbitrary positive
integers by applying the standard inductive
inequality to the largest prime-power divisor,
with the contribution of the remaining factor
becoming negligible after normalization.
Finally, the rank-three result follows by
combining our general transfer theorem with
Zakharov's bound for the weak EGZ constant.

	\begin{proof}[Proof of Theorem~\ref{thm:transfer}]
		Let $W\ge1$ be an integer such that $w(\Z_p^r)\le W$ for all sufficiently large primes $p$. We first prove
		\begin{equation}\label{eq:prime-powers}
			\mathsf s(\Z_N^r)\le(W+o_r(1))N
		\end{equation}
		uniformly as $N\to\infty$ through prime powers, including powers of a fixed small prime.
		
		Assume to the contrary that there are $\eta>0$, prime powers $N_j\to\infty$, and sequences $S_j$ over $\Z_{N_j}^r$ of length $\lfloor(W+\eta)N_j\rfloor$ having no $N_j$-term zero sum.  Set $\delta=\eta/[2(W+\eta)]$.
By Proposition~\ref{prop:decomposition}, there exist
an infinite set $J\subseteq\mathbb N$ and a fixed
finite nonempty index set $E$ such that, for every
$j\in J$, the sequence $S_j$ admits the block
decomposition described in that proposition.
In particular, each block family is complete in
the interior, and, as $j\to\infty$ through $J$,
\[
q_{i,j}=\frac{N_j}{M_{i,j}}\longrightarrow\infty
\qquad(i\in E),
\]
while the limiting normalized lengths $c_i>0$
satisfy
\[
\sum_{i\in E}c_i
\ge(1-\delta)(W+\eta)
=W+\frac{\eta}{2}>W.
\]
Since $N_j$ is a prime power, its divisors $M_{i,j}$
are totally ordered by divisibility. Hence
\[
Q_j=\frac{N_j}{\operatorname{lcm}_{i\in E}M_{i,j}}
=\min_{i\in E}q_{i,j}\longrightarrow\infty,
\]
where the convergence follows because $E$ is fixed
and finite. Thus all the hypotheses of
Proposition~\ref{prop:recombination} are satisfied.
It follows that, for all sufficiently large $j$,
the original sequence $S_j$ contains a zero-sum
subsequence of length $N_j$, contrary to our
assumption. This proves \eqref{eq:prime-powers}.

		For $n>1$, let $Q(n)$ be its largest prime-power divisor and write $n=Q(n)m$. One has $Q(n)\to\infty$ as $n\to\infty$: for fixed $B$, the integers with $Q(n)\le B$ all divide the fixed integer
		\[
		\prod_{\substack{p\le B\\p\ \mathrm{prime}}}
		p^{\lfloor\log_p B\rfloor}.
		\]
		By Lemma~\ref{lem:multiplicative} and \eqref{eq:linear-egz-bound},
		\begin{align*}
			\frac{\mathsf s(\Z_n^r)}n
			&\le\frac{\mathsf s(\Z_{Q(n)}^r)-1}{Q(n)}
			+\frac{\mathsf s(\Z_m^r)}{Q(n)m}\\
			&\le\frac{\mathsf s(\Z_{Q(n)}^r)-1}{Q(n)}
			+\frac{C_r}{Q(n)}
			\le W+o_r(1).
		\end{align*}
		
		Taking $W=\binom{2r-1}{r}+1$, together with Lemma \ref{prop:weak-nine}, we obtain the desired result.
	\end{proof}
	
	\begin{proof}[Proof of Theorem~\ref{thm:main}]
		By Lemma \ref{prop:weak-nine}, $w(\F_p^3)=9$ for all sufficiently large primes $p$. Applying Theorem~\ref{thm:transfer} with $r=3$ and $W=9$ gives $\mathsf s(\Z_n^3)\le9n+o(n)$ for all positive integers $n$. For odd $n$, the well known lower bound $\mathsf s(\Z_n^3)\ge9n-8$ of Elsholtz \cite{Elsholtz} provides the matching lower estimate.
	\end{proof}

\medskip
\noindent\textbf{Statement of AI.} During the preparation of this work,  ChatGPT provided useful simplifications and suggestions for certain proofs, and polished some statements. The authors are completely responsible for the content of the paper.

\medskip
\noindent\textbf{Acknowledgment.} This work was supported by the National Natural Science Foundation of China Grant No.~12571578 and Guangdong
Basic and Applied Basic Research Foundation Grant No.~2024A1515012564.

\end{document}